\documentclass[10pt]{extarticle}

\usepackage[english]{babel}
\usepackage{graphicx}
\usepackage{framed}
\usepackage[normalem]{ulem}
\usepackage{amsmath}
\usepackage{amsthm}
\usepackage{dsfont}
\usepackage{amssymb}
\usepackage{amsfonts}
\usepackage{enumerate}
\usepackage[utf8]{inputenc}
\usepackage[top=1 in,bottom=1 in, left=1 in, right=1 in]{geometry}
\usepackage{mathrsfs}
\usepackage[nottoc, notlof, notlot]{tocbibind}
\usepackage{dsfont}
\usepackage{xcolor}
\usepackage{hyperref}
\usepackage{mathabx}
\usepackage{cases}

\newtheorem{theorem}{Theorem}

\theoremstyle{remark}

\theoremstyle{definition}

\title{Stochastic production planning with regime switching }
\author{
  Elena Cristina Canepa\footnote{Department of Mathematical Methods and Models, University Politehnica of
Bucharest, Romania}\\
  \texttt{cristinacanepa@yahoo.com}
  \and
  Dragos-Patru Covei\footnote{Department of Applied Mathematics, The Bucharest University of
Economic Studies, Piata Romana, 1st district, postal code: 010374, postal
office: 22, Romania}\\
  \texttt{patrucovei@yahoo.com}
	\and
   Traian A. Pirvu\footnote{Department of Mathematics and Statistics, McMaster University, 1280 Main
Street West, Hamilton, ON, L8S 4K1, Canada}\\
  \texttt{tpirvu@math.mcmaster.ca}
}

\begin{document}

\maketitle

\begin{abstract}
This paper considers a stochastic production planning problem with regime
switching. There are two regimes corresponding to different economic cycles.
A factory is planning its production so as to minimize production costs. We
analyze this problem through the value function approach. The optimal
production is characterized through the solution of an elliptic system of
partial differential equations which is shown to have a solution.
\end{abstract}

\section{Introduction}

The purpose of this paper is to consider a stochastic production planning
problem with regime switching parameters and to provide a mathematical
treatment for it. Regime switching modelling is present in many areas such
as financial economics and management. In finance we point the interested
reader to \cite{PZ0}, \cite{PZ} and the references therein.

In the last decade an extensive literature on production planing/management
with regime switching emerged. We only recall a few works. The paper \cite%
{CLP} studies the cost minimization problem of a company within an economy
characterized by two regimes. In civil engineering \cite{DMDK} studies the
optimal stochastic control problem for home energy systems with regime
switching; the two regimes are the peak and off peak energy demand. The work 
\cite{GK} considers the production control problem in a manufacturing system
with multiple machines which are subject to breakdowns and repairs. The
mathematical modelling for these problems makes it possible to find
solutions by simply solving stochastic control problems with regime
dependent controls/value functions. The paper \cite{CP} provides the
mathematical analysis and results of a fairly general class of stochastic
control problems such as the ones appearing in stochastic production
planning over infinite horizons and with regime dependent model parameters.
Their solution approach relies on the concept of value function and the
later is characterized through a system of elliptic equations which is shown
to have solutions. Among recent papers which contribute to the mathematical
analysis of stochastic planning problem we mention \cite{DCT}, \cite{CDPAMC}
and \cite{DT}.

In this paper we look at production planning problem with regime switching
parameters in a random environment. A factory is planning its production of
several economic goods as to minimize inter temporally its production and
inventory costs. A constant discount rate is used to measure on the same
time scale costs which occur at different times. The stochasticity is driven
by a $N$-dimensional Brownian motion and a Markov chain. The Markov chain
models the different economic regimes while the multidimensional Brownian
motion captures the random nature of good's demand; the demand is also
linked to economic cycles and this makes it dependent on the Markov chain as
well. The constant discount rate may also depend on the Markov chain. We add
a stopping criterion in evaluating the inter temporal costs, which is the
stopping time when the inventory of the goods exceeds some threshold level.
We tackle this production planning problem by the value function approach.
Using probabilistic techniques we derive the Hamilton Jacobi Bellman (HJB)
of the value function. We employ partial differential equations (PDE)
tools/techniques to analyze the HJB equation. In the end we prove a
verification result, i.e., we show that the HJB equation yields the optimal
production.

The remainder of this paper is organized as follows. Section \ref{sec2}
presents the model and the objectives. In Section \ref{sec3} we present the
methodology.

\section{Formulation of the model \label{sec2}}

We begin our presentation of the problem to be studied by considering a
factory producing $N$ types of economic goods which are stored in an
inventory designated place. The factory would like to tune its production of
the goods in such a way as to minimize production costs and inventory costs.
We allow for regime switching in our model; regime switching refers to the
situations when the characteristics of the state process are affected by
several regimes (e.g. in finance bull and bear market with higher volatility
in the bear market), or economic cycles characterized by high versus low
demand of economic products (e.g. in the auto industry there is a higher
demand for cars in summer time).

Next, we formulate the model mathematically. There exists a complete
probability space 
\begin{equation*}
(\Omega ,\mathcal{F},\{\mathcal{F}_{t}\}_{0\leq t\leq \infty },P),
\end{equation*}%
on which lives a $N$-dimensional Brownian motion denoted by 
\begin{equation*}
w=\left( w_{1},...,w_{N}\right) .
\end{equation*}%
The regime switching is captured by a continuous time homogeneous Markov
chain $\epsilon (t)$ adapted to $\mathcal{F}_{t}$ with two regimes good and
bad, i.e., $\epsilon (t)\in \{{1},{2}\},\quad t\in \lbrack 0,\infty )$. In a
specific application, $\epsilon (t)=1$ could represent a regime of economic
growth while $\epsilon (t)=2$ could represent a regime of economic
recession. In another application, $\epsilon (t)=1$ could represent a regime
in which consumer demand is high while $\epsilon (t)=2$ could represent a
regime in which consumer demand is low.

The Markov chain's rate matrix is 
\begin{equation}
A=\left( 
\begin{array}{cc}
-a_{1} & a_{1} \\ 
a_{2} & -a_{2}%
\end{array}%
\right) ,  \label{4}
\end{equation}%
for some $a_{1}>0,$ $a_{2}>0$. Diagonal elements $A_{ii}$ are defined such
that%
\begin{equation}
A_{ii}=-\underset{j\neq i}{\Sigma }A_{ij},  \label{5}
\end{equation}%
where 
\begin{equation*}
A_{11}=-a_{1},A_{12}=a_{1},A_{21}=a_{2},A_{22}=-a_{2}.
\end{equation*}%
In this case, if $p_{t}=\mathbb{E}[\epsilon (t)]\in \mathbb{R}^{2}$, then 
\begin{equation}
\frac{d\epsilon (t)}{dt}=A\epsilon (t).  \label{6}
\end{equation}%
Moreover 
\begin{equation}
\epsilon (t)=\epsilon (0)+\int_{0}^{t}A\epsilon (u)\,du+M({t}),  \label{mj}
\end{equation}%
where ${M(t)}$ is a martingale with respect to $\mathcal{F}_{t}$. The
filtration $\mathcal{F}_{t}$ is generated by the $N$-dimensional Brownian
motion and the Markov chain.

Next, let us introduce the control variables in our model. Let 
\begin{equation*}
p\left( t\right) =\left( p_{1}(t, \epsilon (t)),...,p_{N}(t, \epsilon
(t))\right) ,
\end{equation*}%
represent the production rate at time $t$ (control variable) \textbf{%
adjusted for the demand rate}. That means we subtract the demand rate so
that we obtain net production rate. Next, let $y_{i}^{0, \epsilon (0)}$
denote the initial inventory level of good $i,$ and $y_{i}(t, \epsilon (t))$
the inventory level of good $i$, at time $t,$ \textbf{adjusted for demand}.
Again, we look at the net inventory since it is this quantity which incurs
inventory costs. These adjusted for demand inventory levels are modelled by
the following system of stochastic differential equations

\begin{equation}
dy_{i}\left( t,\epsilon (t)\right) =p_{i}dt+\sigma _{\epsilon (t)}dw_{i}%
\text{, }y_{i}\left( 0,\epsilon (0)\right) =y_{i}^{0,\epsilon (0)}\text{, }%
i=1,...,N,  \label{cons0}
\end{equation}%
where $\sigma _{\epsilon (t)}$ is a regime dependent constant (non-zero)
diffusion coefficient taking on two values, $\sigma _{1}$ and $\sigma _{2}$.
The stochasticity here is due to demand adjustment which is random in nature
and dependent on the regime. Another source of randomness our model ca
accommodate are inventory spoilages. One can think of examples when the
demand is more volatile in some periods (e.g. some states of the Markov
chain) and less volatile in other periods.

We impose a stopping production criterion; that is when the (net) inventory
exceeds an exogenous threshold level then the production stops (this is
often the case in auto industry when the storage capacity of newly produced
cars is exhausted). Let us formalize this mathematically; $\tau $ denotes
the stopping time representing the moment when the (net) inventory level
reaches some positive threshold $R$, i.e., 
\begin{equation*}
\tau =\inf_{t>0}\{\left\vert y(t,\epsilon (t))\right\vert \geq R\}.
\end{equation*}%
Here, $|\cdot |$ stands for the Euclidian norm. At this point we are ready
to state our objective.

\subsection{The Objective}

The performance over time of a demand adjusted production rate(s) 
\begin{equation*}
p\left( t,\epsilon (t)\right) =\left( p_{1}(t,\epsilon
(t)),...,p_{N}(t,\epsilon (t))\right)
\end{equation*}%
is measured by means of its production costs and inventory costs. At this
point we introduce the cost functional which measures the quadratic loss: 
\begin{equation}
J\left( p_{1},...,p_{N}\right) :=\text{ }E\int_{0}^{\tau }(|p(t,\epsilon
(t))|^{2}+f_{\epsilon (t)}\left( y(t,\epsilon (t))\right) )e^{-\alpha
_{\epsilon (t)}t}dt,  \label{opti}
\end{equation}%
where $|p(t,\epsilon (t))|^{2}$ and $f_{\epsilon (t)}\left( y(t,\epsilon
(t))\right) $ denote the quadratic holding cost and the production cost
functions, respectively. Again let us recall that we measure deviations from
the demand, whence the loss. Here $\alpha _{\epsilon (t)}$ is a regime
dependent (taking on two values $\alpha _{1}$ and $\alpha _{2}$), constant
psychological rate of time discount, whence the exponential discounting. The
constant psychological rate of time discount is employed to measure on the
same time scale outcomes which occur at different times.

At this point we are ready to frame our objective, which is to minimize the
cost functional. i.e., 
\begin{equation}
\underset{( p_{1},...,p_{N} ) \in \mathbb{R}^{N}}{\inf }\{J\left(
p_{1},...,p_{N}\right) \}\text{, }  \label{min}
\end{equation}%
subject to the stochastic differential equation system (\ref{cons0}).

\section{The Methodology \label{sec3}}

Having presented the problem we want to solve, now we provide our means to
tackle it. Our approach is based on the value function and dynamic
programming which leads to an HJB system of equations.

We apply probabilistic techniques to characterize the value function; that
is we search for functions $z_{i}:R\rightarrow R,\,\,i=1,2$ such that the
stochastic process $Z^{p}(t)$ defined below 
\begin{equation}
Z^{p}\left( t\right) =-e^{-\alpha _{\epsilon (t)}t}z_{\epsilon (t)}\left(
y\left( t,\epsilon (t)\right) \right) -\int_{0}^{t}[|p(s,\epsilon
(s))|^{2}+f_{\epsilon (t)}\left( y(s,\epsilon (s))\right) ]e^{-\alpha
_{\epsilon (s)}s}\,ds,  \label{martsub}
\end{equation}%
is supermartingale for all 
\begin{equation*}
p\left( t,\epsilon (t)\right) =\left( p_{1}(t,\epsilon
(t)),...,p_{N}(t,\epsilon (t))\right) 
\end{equation*}%
and martingale for the optimal control 
\begin{equation*}
p^{\ast }\left( t,\epsilon (t)\right) =\left( {p}_{1}^{\ast }(t,\epsilon
(t)),...,{p}_{N}^{\ast }(t,\epsilon (t))\right) .
\end{equation*}%
Let $B_{R}=\left\{ x\in \mathbb{R}^{N}\left\vert \left\vert x\right\vert
<R\right. \right\} $ be the open ball of radius $R>0$ centered at the
origin. We search for $z_{1}$, $z_{2}$ functions in $C^{2}\left(
B_{R}\right) \cap C\left( \overline{B}_{R}\right) $, and the
supermartingale/martingale requirement yields by means of It\^{o}'s Lemma
for Markov modulated diffusions the HJB system of equations which
characterizes the value function 
\begin{equation}
-a_{1}z_{2}+(a_{1}+\alpha _{1})z_{1}-\frac{\sigma _{1}^{2}}{2}\Delta
z_{1}-f_{1}\left( x\right) =\underset{p\in \mathbb{R}^{N}}{\inf }\{p\nabla
z_{1}+\left\vert p\right\vert ^{2}\text{ }\},  \label{solv1}
\end{equation}%
and 
\begin{equation}
-a_{2}z_{1}+(a_{2}+\alpha _{2})z_{2}-\frac{\sigma _{2}^{2}}{2}\Delta
z_{2}-f_{2}\left( x\right) =\underset{p\in \mathbb{R}^{N}}{\inf }\{p\nabla
z_{2}+\left\vert p\right\vert ^{2}\},  \label{solv2}
\end{equation}%
where $x\in \mathbb{R}^{N}$ assumes values $(y_{1}\left( 0,\varepsilon
(0)\right) ,...,y_{N}\left( 0,\varepsilon (0)\right) )$, $f_{1}$, $f_{2}:%
\overline{B}_{R}\rightarrow \left[ 0,\infty \right) $ are continuous, convex
functions satisfying%
\begin{equation}
\text{there exists }M_{i}>0\text{ such that }f_{i}\left( x\right) \leq
M_{i}\left\vert x\right\vert ^{2}\text{, }i=1,2,  \label{acond}
\end{equation}%
and $\lambda _{\epsilon (t)}$ is a regime dependent (taking on two values $%
\lambda _{1}>0$ and $\lambda _{2}>0$), constant psychological rate of time
discount, whence the exponential discounting.

This HJB system can be turned into a partial differential equation system
(PDE system) since a simple calculation yields 
\begin{equation}
\underset{p\in \mathbb{R}^{N}}{\inf }\{p\nabla z_{j}+\left\vert p\right\vert
^{2}\}=-\frac{1}{4}\left\vert \nabla z_{j}\right\vert ^{2},\,\,j=1,2.
\label{foc}
\end{equation}%
Thus, the HJB system becomes the PDE system 
\begin{equation}
\left\{ 
\begin{array}{c}
-a_{1}z_{2}+(a_{1}+\alpha _{1})z_{1}-\frac{{\sigma _{1}}^{2}}{2}\Delta
z_{1}-f_{1}\left( x\right) =-\frac{1}{4}\left\vert \nabla z_{1}\right\vert
^{2}\text{ for }x\in B_{R}\text{,} \\ 
-a_{2}z_{1}+(a_{2}+\alpha _{2})z_{2}-\frac{{\sigma _{2}}^{2}}{2}\Delta
z_{2}-f_{2}\left( x\right) =-\frac{1}{4}\left\vert \nabla z_{2}\right\vert
^{2}\text{ for }x\in B_{R}\text{.}%
\end{array}%
\right.  \label{sysd}
\end{equation}%
In order to perform the verification, i.e., show that the HJB system gives
the solution of the optimization problem, one needs to impose the following
boundary condition 
\begin{equation}
z_{1}(x)=z_{2}(x)=0\text{ for }x\in \partial B_{R}.  \label{bdc}
\end{equation}%
The gradient term in the above PDE system can be removed by the change of
variable 
\begin{equation*}
u_{j}\left( x\right) =e^{\frac{-z_{j}\left( x\right) }{2\sigma _{j}^{2}}%
},\quad j=1,2,
\end{equation*}%
to get a simpler PDE system

\begin{equation}
\left\{ 
\begin{array}{l}
\Delta u_{1}\left( x\right) =u_{1}\left( x\right) [\frac{1}{{\sigma _{1}^{4}}%
}f_{1}\left( x\right) +\frac{2(a_{1}+\alpha _{1})}{{\sigma _{1}^{2}}}\ln
u_{1}\left( x\right) -2a_{1}\frac{{\sigma _{2}^{2}}}{{\sigma _{1}^{4}}}\ln
u_{2}\left( x\right) ]\text{ for }x\in B_{R}, \\ 
\Delta u_{2}\left( x\right) =u_{2}\left( x\right) [\frac{1}{{\sigma }_{2}^{4}%
}f_{2}\left( x\right) +\frac{2(a_{2}+\alpha _{2})}{{\sigma _{2}^{2}}}\ln
u_{2}\left( x\right) -2a_{2}\frac{{\sigma _{1}^{2}}}{{\sigma _{2}^{4}}}\ln
u_{1}\left( x\right) ]\text{ for }x\in B_{R}\text{,} \\ 
u_{2}\left( x\right) >0\text{, }u_{1}\left( x\right) >0\text{ for }x\in B_{R}%
\text{,}%
\end{array}%
\right.  \label{I1}
\end{equation}%
with the corresponding boundary condition 
\begin{equation}
u_{1}(x)=u_{2}(x)=1\text{ for }x\in \partial B_{R}.  \label{bdc1}
\end{equation}%
The value function will give us in turn the candidate optimal control. The
first order optimality conditions on the lefthand side of (\ref{foc}) are
sufficient for optimality since we deal with a quadratic function to be
optimized, and they produce the candidate optimal control as follows:

\begin{equation*}
{p}_{i}^{\ast }(t,\epsilon (t))=\overline{p}_{i}(y_{1}\left( t,\epsilon
(t)\right) ,\ldots ,y_{N}\left( t,\epsilon (t)\right) )\text{, }i=1,...,N,
\end{equation*}%
and 
\begin{equation}
\overline{p}_{i}(x_{1},...,x_{N},j)=-\frac{1}{2}\frac{\partial z_{j}}{%
\partial x_{i}}\left( x_{1},...,x_{N}\right) \text{, for }i=1,...,n,\,\,j=1,2%
\text{.}  \label{optio}
\end{equation}

The system (\ref{sysd})-(\ref{bdc}) is key in solving our problem so we need
to analyze it. We prove the following result:

\begin{theorem}
\label{sis}The system of equations (\ref{sysd})-(\ref{bdc}) has a unique
positive classical solution $\left( z_{1},z_{2}\right) $. Moreover, 
\begin{equation}
z_{i}\left( x\right) \leq -2\sigma _{i}^{2}K_{i}(R^{2}-|x|^{2})\text{, for
some }K_{i}<0,\quad i=1,2.  \label{sq1}
\end{equation}
\end{theorem}

\subparagraph{\textbf{Proof }}

Our approach, being constructive, will be useful for a computational scheme
for numerical approximations of the solution. Since the system (\ref{sysd})-(%
\ref{bdc}) is equivalent to (\ref{I1})-(\ref{bdc1}) we will work with the
later. We proceed in three steps: step 1) establishes a sub-solution and a
super-solution; step 2) provides an approximating sequence of functions
which converges to the solution; step 3) established the uniqueness of the
solution.

\textbf{Step 1} The main problem is reduced to the construction of the
function $\left( \underline{u}_{1},\underline{u}_{2}\right) $ called
sub-solution and a function $\left( \overline{u}_{1},\overline{u}_{2}\right) 
$ named super-solution with order (i.e., $\underline{u}_{1}\left( x\right)
\leq \overline{u}_{1}\left( x\right) $ and $\underline{u}_{2}\left( x\right)
\leq \overline{u}_{2}\left( x\right) $, for all $x\in \overline{B}_{R}$ ) to
the system (\ref{I1}), which satisfy the inequalities%
\begin{equation}
\left\{ 
\begin{array}{l}
\Delta \underline{u}_{1}\left( x\right) \geq \underline{u}_{1}\left(
x\right) [\frac{1}{{\sigma _{1}^{4}}}f_{1}\left( x\right) +\frac{%
2(a_{1}+\alpha _{1})}{{\sigma _{1}^{2}}}\ln \underline{u}_{1}\left( x\right)
-\frac{2a_{1}{\sigma _{2}^{2}}}{{\sigma _{1}^{4}}}\ln \underline{u}%
_{2}\left( x\right) ]\text{, }x\in B_{R}, \\ 
\Delta \underline{u}_{2}\left( x\right) \geq \underline{u}_{2}\left(
x\right) [\frac{1}{{\sigma _{2}}^{4}}f_{2}\left( x\right) +\frac{%
2(a_{2}+\alpha _{2})}{{\sigma _{2}^{2}}}\ln \underline{u}_{2}\left( x\right)
-\frac{2a_{2}{\sigma _{1}^{2}}}{{\sigma _{2}^{4}}}\ln \underline{u}%
_{1}\left( x\right) ]\text{, }x\in B_{R}, \\ 
\Delta \overline{u}_{1}\left( x\right) \leq \overline{u}_{1}\left( x\right) [%
\frac{1}{{\sigma _{1}^{4}}}f_{1}\left( x\right) +\frac{2(a_{1}+\alpha _{1})}{%
{\sigma _{1}^{2}}}\ln \overline{u}_{1}\left( x\right) -\frac{2a_{1}{\sigma
_{2}^{2}}}{{\sigma _{1}^{4}}}\ln \overline{u}_{2}\left( x\right) ]\text{, }%
x\in B_{R}, \\ 
\Delta \overline{u}_{2}\left( x\right) \leq \overline{u}_{2}\left( x\right) [%
\frac{1}{{\sigma _{2}^{4}}}f_{2}\left( x\right) +\frac{2(a_{2}+\alpha _{2})}{%
{\sigma _{2}^{2}}}\ln \overline{u}_{2}\left( x\right) -\frac{2a_{2}{\sigma
_{1}^{2}}}{{\sigma _{2}^{4}}}\ln \overline{u}_{1}\left( x\right) ]\text{, }%
x\in B_{R}.%
\end{array}%
\right.  \label{sissu}
\end{equation}%
The construction of the sub-solution requires some work. More exactly, by
direct calculations we observe that there exist 
\begin{equation}
\left( \underline{u}_{1}\left( x\right) ,\underline{u}_{2}\left( x\right)
\right) =\left( e^{K_{1}\left( R^{2}-\left\vert x\right\vert ^{2}\right)
},e^{K_{2}\left( R^{2}-\left\vert x\right\vert ^{2}\right) }\right) ,\text{
with }K_{1},K_{2}\in \left( -\infty ,0\right) ,  \label{subs}
\end{equation}%
satisfying (\ref{sissu}). By substituting (\ref{subs}) into (\ref{sissu}) we
prove that there exist $K_{1},K_{2}\in \left( -\infty ,0\right) $ such that

\begin{equation*}
\left\{ 
\begin{array}{l}
4K_{1}^{2}\left\vert x\right\vert ^{2}-2K_{1}N\allowbreak \geq \frac{M_{1}}{{%
\sigma _{1}^{4}}}\left\vert x\right\vert ^{2}+\frac{2(a_{1}+\alpha _{1})K_{1}%
}{{\sigma _{1}^{2}}}\left( R^{2}-\left\vert x\right\vert ^{2}\right) -\frac{%
2a_{1}{\sigma _{2}^{2}}K_{2}}{{\sigma _{1}^{4}}}\left( R^{2}-\left\vert
x\right\vert ^{2}\right) , \\ 
4K_{2}^{2}\left\vert x\right\vert ^{2}-2K_{2}N\allowbreak \geq \frac{M_{2}}{{%
\sigma _{2}^{4}}}\left\vert x\right\vert ^{2}+\frac{2(a_{2}+\alpha _{2})K_{2}%
}{{\sigma _{2}^{2}}}\left( R^{2}-\left\vert x\right\vert ^{2}\right) -\frac{%
2a_{2}{\sigma _{1}^{2}}K_{1}}{{\sigma _{2}^{4}}}\left( R^{2}-\left\vert
x\right\vert ^{2}\right) ,%
\end{array}%
\right.
\end{equation*}%
or, equivalently%
\begin{equation*}
\left\{ 
\begin{array}{c}
\lbrack 4K_{1}^{2}-\frac{2a_{1}{\sigma _{2}^{2}}K_{2}+M_{1}-2(a_{1}+\alpha
_{1})K_{1}{\sigma _{1}^{2}}}{{\sigma _{1}^{4}}}]\left\vert x\right\vert ^{2}-%
\frac{2(a_{1}+\alpha _{1})K_{1}}{{\sigma _{1}^{2}}}R^{2}-2K_{1}N+\frac{2a_{1}%
{\sigma _{2}^{2}}K_{2}}{{\sigma _{1}^{4}}}R^{2}\geq 0, \\ 
\lbrack 4K_{2}^{2}-\frac{2a_{2}{\sigma _{1}^{2}}K_{1}+M_{2}-2(a_{2}+\alpha
_{2})K_{2}{\sigma _{2}^{2}}}{{\sigma _{2}^{4}}}]\left\vert x\right\vert ^{2}-%
\frac{2(a_{2}+\alpha _{2})K_{2}}{{\sigma _{2}^{2}}}R^{2}-2K_{2}N+\frac{2a_{2}%
{\sigma _{1}^{2}}K_{1}}{{\sigma _{2}^{4}}}R^{2}\geq 0.%
\end{array}%
\right.
\end{equation*}%
Therefore, it suffices to show that there exist $K_{1},K_{2}\in \left(
-\infty ,0\right) $ such that 
\begin{equation}
\left\{ 
\begin{array}{l}
4K_{1}^{2}+\frac{2(a_{1}+\alpha _{1}){\sigma _{1}^{2}}}{{\sigma _{1}^{4}}}%
K_{1}-\frac{M_{1}}{{\sigma _{1}^{4}}}-\frac{2a_{1}{\sigma _{2}^{2}}}{{\sigma
_{1}^{4}}}K_{2}\geq 0 \\ 
4K_{2}^{2}+\frac{2(a_{2}+\alpha _{2}){\sigma _{2}^{2}}}{{\sigma _{2}^{4}}}%
K_{2}-\frac{M_{2}}{{\sigma _{2}^{4}}}-\frac{2a_{2}{\sigma _{1}^{2}}}{{\sigma
_{2}^{4}}}K_{1}\geq 0 \\ 
-\frac{2(a_{1}+\alpha _{1})R^{2}}{{\sigma _{1}^{2}}}K_{1}-2K_{1}N+\frac{%
2a_{1}{\sigma _{2}^{2}}R^{2}}{{\sigma _{1}^{4}}}K_{2}\geq 0 \\ 
-\frac{2(a_{2}+\alpha _{2})R^{2}}{{\sigma _{2}^{2}}}K_{2}-2NK_{2}+\frac{%
2a_{2}{\sigma _{1}^{2}}R^{2}}{{\sigma _{2}^{4}}}K_{1}\geq 0.%
\end{array}%
\right.  \label{ineq}
\end{equation}%
To do this, first we take$\allowbreak $%
\begin{equation}
K_{2}=-\max \{\frac{1}{4\sigma _{2}^{2}}(\alpha _{2}+a_{2}+\sqrt{\left(
\alpha _{2}+a_{2}\right) ^{2}+4M_{2}}),\frac{\frac{1}{4\sigma _{1}^{2}}%
(\alpha _{1}+a_{1}+\sqrt{\left( \alpha _{1}+a_{1}\right) ^{2}+4M_{1}})}{2%
\frac{R^{2}\sigma _{2}^{2}K_{2}a_{1}}{\sigma _{1}^{4}}[\frac{2(a_{1}+\alpha
_{1})}{{\sigma _{1}^{2}}}R^{2}+2N]^{-1}}\}  \label{k2}
\end{equation}%
and second 
\begin{equation}
-K_{1}\in \left[ \frac{-2\frac{R^{2}}{\sigma _{1}^{4}}\sigma
_{2}^{2}K_{2}a_{1}}{\frac{2(a_{1}+\alpha _{1})}{{\sigma _{1}^{2}}}R^{2}+2N},%
\frac{-[\frac{2(a_{2}+\alpha _{2})}{{\sigma _{2}^{2}}}R^{2}+2N]K_{2}}{2\frac{%
R^{2}}{\sigma _{2}^{4}}\sigma _{1}^{2}a_{2}}\right] .  \label{k1}
\end{equation}%
Then, taking into account (\ref{k2}) and (\ref{k1}), we see that 
\begin{eqnarray*}
\left[ \frac{2(a_{2}+\alpha _{2})}{{\sigma _{2}^{2}}}R^{2}+2N\right] \left[ 
\frac{2(a_{1}+\alpha _{1})}{{\sigma _{1}^{2}}}R^{2}+2N\right] &\geq &\frac{%
2(a_{2}+\alpha _{2})}{{\sigma _{2}^{2}}}R^{2}\frac{2(a_{1}+\alpha _{1})}{{%
\sigma _{1}^{2}}}R^{2} \\
&\geq &\frac{4a_{1}a_{2}}{{\sigma _{1}^{2}\sigma _{2}^{2}}}R^{4} \\
&=&\left( 2a_{1}\frac{{\sigma _{2}^{2}}}{{\sigma _{1}^{4}}}R^{2}\right)
\left( 2a_{2}\frac{{\sigma _{1}^{2}}}{{\sigma _{2}^{4}}}R^{2}\right)
\end{eqnarray*}%
and, consequently, we obtain%
\begin{equation}
\frac{-2a_{1}\frac{{\sigma _{2}^{2}}}{{\sigma _{1}^{4}}}R^{2}K_{2}}{\frac{%
2(a_{1}+\alpha _{1})}{{\sigma _{1}^{2}}}R^{2}+2N}\leq \frac{-\left[ \frac{%
2(a_{2}+\alpha _{2})}{{\sigma _{2}^{2}}}R^{2}+2N\right] K_{2}}{2a_{2}\frac{{%
\sigma _{1}^{2}}}{{\sigma _{2}^{4}}}R^{2}}.  \label{ink1}
\end{equation}%
The inequality (\ref{ink1}) says that it is possible to choose $K_{1}$ as in
(\ref{k1}). On the other hand, in (\ref{k2}) and (\ref{k1}), it follows that%
\begin{equation*}
\frac{1}{4\sigma _{1}^{2}}[\alpha _{1}+a_{1}+\sqrt{\left( \alpha
_{1}+a_{1}\right) ^{2}+4M_{2}}]\leq \frac{-2a_{1}\frac{{\sigma _{2}^{2}}}{{%
\sigma _{1}^{4}}}R^{2}K_{2}}{\frac{2(a_{1}+\alpha _{1})}{{\sigma _{1}^{2}}}%
R^{2}+2N}\leq -K_{1}.
\end{equation*}%
One only has to notice that:

1.%
\begin{equation*}
4K_{1}^{2}+\frac{2(a_{1}+\alpha _{1})}{{\sigma _{1}^{2}}}K_{1}-\frac{M_{1}}{{%
\sigma _{1}^{4}}}\geq 0
\end{equation*}%
since%
\begin{equation*}
\frac{1}{4\sigma _{1}^{2}}[\alpha _{1}+a_{1}+\sqrt{\left( \alpha
_{1}+a_{1}\right) ^{2}+4M_{1}}]\leq -K_{1};
\end{equation*}

2.%
\begin{equation*}
4K_{2}^{2}+\frac{2(a_{2}+\alpha _{2})}{{\sigma _{2}^{2}}}K_{2}-\frac{M_{2}}{{%
\sigma _{2}^{4}}}\geq 0
\end{equation*}%
since%
\begin{equation*}
\frac{1}{4\sigma _{2}^{2}}[\alpha _{2}+a_{2}+\sqrt{\left( \alpha
_{2}+a_{2}\right) ^{2}+4M_{2}}]\leq -K_{2};
\end{equation*}

3.%
\begin{equation*}
-[\frac{2(a_{1}+\alpha _{1})}{{\sigma _{1}^{2}}}R^{2}+2N]K_{1}+2a_{1}\frac{{%
\sigma _{2}^{2}}}{{\sigma _{1}^{4}}}R^{2}K_{2}\geq 0
\end{equation*}%
since%
\begin{equation*}
-K_{1}\geq \frac{-2a_{1}\frac{{\sigma _{2}^{2}}}{{\sigma _{1}^{4}}}R^{2}K_{2}%
}{\frac{2(a_{1}+\alpha _{1})}{{\sigma _{1}^{2}}}R^{2}+2N};
\end{equation*}

4.%
\begin{equation*}
-[\frac{2(a_{2}+\alpha _{2})}{{\sigma _{2}^{2}}}R^{2}+2N]K_{2}+2a_{2}\frac{{%
\sigma _{1}^{2}}}{{\sigma _{2}^{4}}}R^{2}K_{1}\geq 0
\end{equation*}%
since%
\begin{equation*}
-K_{1}\leq \frac{-\left[ \frac{2(a_{2}+\alpha _{2})}{{\sigma _{2}^{2}}}%
R^{2}+2N\right] K_{2}}{2a_{2}\frac{{\sigma _{1}^{2}}}{{\sigma _{2}^{4}}}R^{2}%
}
\end{equation*}
and, thus (\ref{subs}) is a sub-solution for the system (\ref{I1}).
Constructing a super-solution is easier. It turns out that

\begin{equation*}
\text{ }\left( \overline{u}_{1}\left( x\right) ,\text{ }\overline{u}%
_{2}\left( x\right) \right) =\left( 1,1\right) ,
\end{equation*}%
is a super-solution of (\ref{I1}).

\textbf{Step 2} By the above construction one gets 
\begin{equation*}
\underline{u}_{1}\left( x\right) \leq \overline{u}_{1}\left( x\right) \text{
and }\underline{u}_{2}\left( x\right) \leq \overline{u}_{2}\left( x\right) 
\text{ for all }x\in \overline{B}_{R}\text{.}
\end{equation*}%
Next, we are showing that the problem (\ref{I1}) admits a unique solution 
\begin{equation*}
\left( u_{1},u_{2}\right) \in \lbrack C^{2}\left( B_{R}\right) \cap C\left( 
\overline{B}_{R}\right) ]^{2}=[C^{2}\left( B_{R}\right) \cap C\left( 
\overline{B}_{R}\right) ]\times \lbrack C^{2}\left( B_{R}\right) \cap
C\left( \overline{B}_{R}\right) ],
\end{equation*}%
such that%
\begin{equation*}
\underline{u}_{1}\left( x\right) \leq u_{1}\left( x\right) \leq \overline{u}%
_{1}\left( x\right) \text{ and }\underline{u}_{2}\left( x\right) \leq
u_{2}\left( x\right) \leq \overline{u}_{2}\left( x\right) ,\text{ for }x\in 
\overline{B}_{R}.
\end{equation*}%
Denote%
\begin{equation*}
M_{1}=e^{K_{1}R^{2}}\text{, }M_{2}=e^{K_{2}R^{2}}\text{ and }M=1.
\end{equation*}%
Let $g_{1}:\overline{B}_{R}\times \lbrack M_{1},M]\times \lbrack
M_{2},M]\rightarrow \mathbb{R}$ and $g_{2}:\overline{B}_{R}\times \lbrack
M_{2},M]\times \lbrack M_{1},M]\rightarrow \mathbb{R}$ defined by 
\begin{eqnarray*}
g_{1}\left( x,t,s\right) &=&\frac{1}{\sigma ^{4}}f_{1}\left( x\right) t+%
\frac{2(a_{1}+\alpha _{1})}{{\sigma _{1}^{2}}}t\ln t-\frac{2a_{1}{\sigma
_{2}^{2}}}{{\sigma _{1}^{4}}}t\ln s, \\
g_{2}\left( x,t,s\right) &=&\frac{1}{{\sigma }_{2}^{4}}f_{2}\left( x\right)
s+\frac{2(a_{2}+\alpha _{2})}{{\sigma _{2}^{2}}}s\ln s-\frac{2a_{2}{\sigma
_{1}^{2}}}{{\sigma _{2}^{4}}}s\ln t.
\end{eqnarray*}%
Since $g_{1}$ is a continuous function with respect to the first variable in 
$\overline{B}_{R}$ and continuously differentiable with respect to the
second and third in $[M_{1},M]\times \lbrack M_{2},M],$ it allows to choose $%
\Lambda _{1}\in \left( -\infty ,0\right) $ such that 
\begin{equation*}
-\Lambda _{1}\geq \frac{g_{1}\left( x,t_{1},s\right) -g_{1}\left(
x,t_{2},s\right) }{t_{2}-t_{1}}\text{, }
\end{equation*}%
for every $t_{1},t_{2}$ with $\underline{u}_{1}\leq t_{2}<t_{1}\leq 
\overline{u}_{1}$ and $x\in B_{R}$. Similarly for $g_{2}$, one can set $%
\Lambda _{2}\in \left( -\infty ,0\right) $ such that%
\begin{equation*}
-\Lambda _{2}\geq \frac{g_{2}\left( x,t,s_{1}\right) -g_{2}\left(
x,t,s_{2}\right) }{s_{2}-s_{1}},
\end{equation*}%
for every $s_{1},s_{2}$ with $\underline{u}_{2}\leq s_{2}<s_{1}\leq 
\overline{u}_{2}$ and $x\in B_{R}$.

We develop a sequence of approximations for the solution. The sub- and
super- solution will be used as the initial iteration in a Picard type of
monotone iteration process. Namely, with the starting point $\left(
u_{1}^{0},u_{2}^{0}\right) =\left( \underline{u}_{1},\underline{u}%
_{2}\right) $ we inductively define a sequence $\left\{ \left(
u_{1}^{k},u_{2}^{k}\right) \right\} _{k\in \mathbb{N}^{\ast }}$ such that 
\begin{equation*}
\left\{ 
\begin{array}{ll}
\Delta u_{1}^{k}+\Lambda _{1}u_{1}^{k}=g_{1}\left(
x,u_{1}^{k-1},u_{2}^{k-1}\right) +\Lambda _{1}u_{1}^{k-1} & \text{for }x\in
B_{R}, \\ 
\Delta u_{2}^{k}+\Lambda _{2}u_{2}^{k}=g_{2}\left(
x,u_{1}^{k-1},u_{2}^{k-1}\right) +\Lambda _{2}u_{2}^{k-1} & \text{for }x\in
B_{R}, \\ 
u_{1}^{k}\left( x\right) =u_{2}^{k}\left( x\right) =1 & \text{for }x\in
\partial B_{R}.%
\end{array}%
\right.
\end{equation*}%
The existence proof for 
\begin{equation*}
\left\{ 
\begin{array}{ll}
\Delta u_{1}^{1}+\Lambda _{1}u_{1}^{1}=g_{1}\left(
x,u_{1}^{0},u_{2}^{0}\right) +\Lambda _{1}u_{1}^{0} & \text{for }x\in B_{R},
\\ 
\Delta u_{2}^{1}+\Lambda _{2}u_{2}^{1}=g_{2}\left(
x,u_{1}^{0},u_{2}^{0}\right) +\Lambda _{2}u_{2}^{0} & \text{for }x\in B_{R},
\\ 
u_{1}^{1}\left( x\right) =u_{2}^{1}\left( x\right) =1 & \text{for }x\in
\partial B_{R}.%
\end{array}%
\right.
\end{equation*}%
can be found in \cite{GT}. Clearly, the sequence $\left\{ \left(
u_{1}^{k},u_{2}^{k}\right) \right\} _{k\in \mathbb{N}^{\ast }}$ is well
defined. Next, assuming that 
\begin{equation*}
u_{1}^{k-1}\leq u_{1}^{k}\text{ and }u_{2}^{k-1}\leq u_{2}^{k}\text{ on }%
\overline{B}_{R}\text{,}
\end{equation*}%
we prove that 
\begin{equation*}
u_{1}^{k}\leq u_{1}^{k+1}\text{ and }u_{2}^{k}\leq u_{2}^{k+1}\text{ on }%
\overline{B}_{R}.
\end{equation*}%
The constants $\Lambda _{1}$ and $\Lambda _{2}$ are chosen so that%
\begin{equation}
\begin{array}{ll}
\left( \Delta +\Lambda _{1}\right) \left( u_{1}^{k+1}\left( x\right)
-u_{1}^{k}\left( x\right) \right) \leq 0 & \text{in }B_{R}, \\ 
\left( \Delta +\Lambda _{2}\right) \left( u_{2}^{k+1}\left( x\right)
-u_{2}^{k}\left( x\right) \right) \leq 0 & \text{in }B_{R},%
\end{array}
\label{ump}
\end{equation}%
if 
\begin{equation*}
u_{1}^{k-1}\leq u_{1}^{k}\text{ and }u_{2}^{k-1}\leq u_{2}^{k}\text{ on }%
B_{R}\text{,}
\end{equation*}%
for $k=1,2,...$, which is true for $k=1$ and thus by mathematical induction
for every larger $k$ by (\ref{ump}) and the maximum principle.

Consequently, by induction we get a monotone increasing sequence $\left\{
\left( u_{1}^{k},u_{2}^{k}\right) \right\} _{k\in \mathbb{N}^{\ast }}$ of
iterates 
\begin{eqnarray*}
\underline{u}_{1} &\leq &u_{1}^{1}\leq u_{1}^{2}\leq ...\leq u_{1}^{k-1}\leq
u_{1}^{k}\leq u_{1}^{k+1}\leq ...\leq \overline{u}_{1}\text{ on }\overline{B}%
_{R}\text{,} \\
\underline{u}_{2} &\leq &u_{2}^{1}\leq u_{2}^{2}\leq ...\leq u_{2}^{k-1}\leq
u_{2}^{k}\leq u_{2}^{k+1}\leq ...\leq \overline{u}_{2}\text{ on }\overline{B}%
_{R}\text{.}
\end{eqnarray*}%
To sum up, we have constructed a monotonic and bounded sequence \{$\left(
u_{1}^{k},u_{2}^{k}\right) $\}$_{k\in \mathbb{N}}$ that converges%
\begin{equation*}
\lim_{k\rightarrow \infty }\left( u_{1}^{k}\left( x\right) ,u_{2}^{k}\left(
x\right) \right) =\left( u_{1}\left( x\right) ,u_{2}\left( x\right) \right) 
\text{, for all }x\in \overline{B}_{R}.
\end{equation*}%
Clearly, the limit function $\left( u_{1}\left( x\right) ,u_{2}\left(
x\right) \right) $ exists as a continuous function on $\overline{B}_{R}$.
Via standard bootstrap arguments (see \cite[p. 26]{DS}) 
\begin{equation*}
\left( u_{1}^{k},u_{2}^{k}\right) \overset{k\rightarrow \infty }{\rightarrow 
}\left( u_{1},u_{2}\right) \text{ in }[C^{2}\left( B_{R}\right) \cap C(%
\overline{B}_{R})]^{2},
\end{equation*}%
and $\left( u_{1},u_{2}\right) $ is a solution of problem (\ref{I1})
satisfying 
\begin{equation*}
\underline{u}_{1}\left( x\right) \leq u_{1}\left( x\right) \leq \overline{u}%
_{1}\left( x\right) \text{ and }\underline{u}_{2}\left( x\right) \leq
u_{2}\left( x\right) \leq \overline{u}_{2}\left( x\right) \text{ for all }%
x\in \overline{B}_{R}.
\end{equation*}%
Then 
\begin{equation*}
\left( z_{1}\left( x\right) ,z_{2}\left( x\right) \right) =\left( -2\sigma
_{1}\ln u_{1}\left( x\right) ,-2\sigma _{2}\ln u_{2}\left( x\right) \right)
\in \lbrack C^{2}\left( B_{R}\right) \cap C(\overline{B}_{R})]^{2}\text{,}
\end{equation*}%
is the positive solution of (\ref{sysd})-(\ref{bdc1}) with quadratic growth.

\textbf{Step 3} Let $\left( u_{1},u_{2}\right) $ and $\left( \widetilde{u}%
_{1},\widetilde{u}_{2}\right) $ be any classical positive solutions to the
system (\ref{I1})-(\ref{bdc1}). Clearly%
\begin{equation*}
\lim_{x\rightarrow \partial B_{R}}\frac{u_{1}\left( x\right) }{\widetilde{u}%
_{1}\left( x\right) }=1\text{ and }\lim_{x\rightarrow \partial B_{R}}\frac{%
u_{2}\left( x\right) }{\widetilde{u}_{2}\left( x\right) }=1.
\end{equation*}%
Using the maximum principle coupled with the scalar case in \cite{CP1} we
have 
\begin{equation*}
u_{1}\left( x\right) \leq \widetilde{u}_{1}\left( x\right) \text{ and }%
u_{2}\left( x\right) \leq \widetilde{u}_{2}\left( x\right) ,\text{ for any }%
x\in \overline{B}_{R}.
\end{equation*}%
Thus, by interchanging the roles of $\left( u_{1},u_{2}\right) $ and $\left( 
\widetilde{u}_{1},\widetilde{u}_{2}\right) $, we also have 
\begin{equation*}
\widetilde{u}_{1}\left( x\right) \leq u_{1}\left( x\right) \text{ and }%
\widetilde{u}_{2}\left( x\right) \leq u_{2}\left( x\right) ,\text{ for any }%
x\in \overline{B}_{R}.
\end{equation*}%
It follows that 
\begin{equation*}
\left( u_{1}\left( x\right) ,u_{2}\left( x\right) \right) =\left( \widetilde{%
u}_{1}\left( x\right) ,\widetilde{u}_{2}\left( x\right) \right) ,\text{ for
any }x\in \overline{B}_{R},
\end{equation*}%
thus proving the uniqueness of solution for the problem (\ref{I1})-(\ref%
{bdc1}). The references \cite[Section 3]{MM}, \cite{BM} guarantee also the
uniqueness of a classical solution of the system (\ref{sysd})-(\ref{bdc}).
This completes our proof.

\subsection{Verification}

In this subsection we show that the control of \eqref{optio} is indeed
optimal. This is formalized in the following Theorem.

\begin{theorem}
\label{opt} The production rate(s) ${p}_{i}^{\ast }(t,\epsilon (t)), \,\,
i=1,2,\ldots N$ defined in (\ref{optio}) is optimal. That is for every
production rate(s) $\left( p_{1}(t,\epsilon(t)),...,p_{N}(t,\epsilon
(t))\right)$ 
\begin{equation*}
J\left( p_{1},...,p_{N}\right)\geq J\left( p_{1}^{\ast },...,p_{N}^{\ast
}\right).
\end{equation*}
\end{theorem}

\subparagraph{\textbf{Proof }}

Let us denote by $y^{\ast }(t,\epsilon (t))$ the vector of inventory levels
associated with ${p}_{i}^{\ast }(t,\epsilon (t)),\,\,i=1,2,\ldots N$. Recall
that%
\begin{equation*}
\tau ^{\ast }=\inf_{t>0}\{\left\vert y^{\ast }(t,\epsilon (t))\right\vert
\geq R\},
\end{equation*}%
and%
\begin{equation*}
\tau =\inf_{t>0}\{\left\vert y(t,\epsilon (t))\right\vert \geq R\}.
\end{equation*}%
We proceed in two steps: step 1) shows that the stochastic process $Z^{p}(t)$
defined in \eqref{martsub} is supermartingale for all%
\begin{equation*}
p\left( t,\epsilon (t)\right) =\left( p_{1}(t,\epsilon
(t)),...,p_{N}(t,\epsilon (t))\right) ,
\end{equation*}%
on $0\leq t\leq \tau ,$ and martingale for 
\begin{equation*}
p^{\ast }\left( t,\epsilon (t)\right) =\left( {p}_{1}^{\ast }(t,\epsilon
(t)),...,{p}_{N}^{\ast }(t,\epsilon (t))\right) ,
\end{equation*}%
on $0\leq t\leq \tau ^{\ast }.$ In step 2) we establish optimality of ${p}%
_{i}^{\ast }(t,\epsilon (t)),\,\,i=1,2,\ldots N$ defined in (\ref{optio}).

\textbf{Step 1} It\^{o}'s Lemma for Markov modulated diffusion (see \cite%
{YZY} for more on this) yields%
\begin{eqnarray*}
dZ^{p}\left( t\right) &=&-e^{-\alpha _{\epsilon (s)}s}\bigg[\frac{\sigma
_{\epsilon \left( s\right) }^{2}}{2}\Delta z_{\epsilon (s)}\left( y\left(
s,\epsilon (s)\right) \right) +f_{\epsilon (t)}\left( y\left( s,\epsilon
(s)\right) \right) +p\left( s,\epsilon \left( s\right) \right) \nabla
z_{\epsilon (s)}\left( y\left( s,\epsilon (s)\right) \right) + \\
&&+\left\vert \left( p\left( s,\epsilon (s)\right) \right) \right\vert
^{2}-\left( \alpha _{\epsilon (s)}+a_{\epsilon \left( s\right) }\right)
z_{\epsilon (s)}y\left( s,\epsilon (s)\right) +a_{\epsilon (s)}1_{\left\{
i\neq \epsilon \left( s\right) \right\} }z_{i}\left( y\left( s,\epsilon
(s)\right) \right) \bigg]ds \\
&&-e^{-\alpha _{\epsilon (s)}s}\sigma _{\epsilon \left( s\right) }p\left(
s,\epsilon \left( s\right) \right) \nabla z_{\epsilon (s)}\left( y\left(
s,\epsilon (s)\right) \right) dw\left( s\right) .
\end{eqnarray*}%
Then, the claim yields in light of HJB equation (\ref{solv1}) and (\ref%
{solv2}).

\textbf{Step 2} In the second step let us establish the optimality of $%
\left( p_{1}^{\ast },...,p_{N}^{\ast }\right) $. The
martingale/supermartingale principle yields%
\begin{equation*}
Ee^{-\alpha _{\epsilon \left( \tau ^{\ast }\right) }\tau ^{\ast
}}z_{\epsilon \left( \tau ^{\ast }\right) }\left( y^{\ast }\left( \tau
^{\ast },\epsilon \left( \tau ^{\ast }\right) \right) \right)
+E\int_{0}^{\tau \ast }e^{-\alpha _{\epsilon (u)}u}\left[ \left\vert p^{\ast
}\left( u,\epsilon \left( u\right) \right) \right\vert ^{2}+f_{\epsilon
(t)}\left( y^{\ast }\left( u,\epsilon \left( u\right) \right) \right) \right]
du=z_{\epsilon \left( 0\right) }\left( y^{\ast }\left( 0,\epsilon \left(
0\right) \right) \right) ,
\end{equation*}%
and 
\begin{equation*}
Ee^{-\alpha _{\epsilon \left( \tau \right) }\tau }z_{\epsilon \left( \tau
\right) }\left( y\left( \tau ,\epsilon \left( \tau \right) \right) \right)
+E\int_{0}^{\tau }e^{-\alpha _{\epsilon (u)}u}\left[ \left\vert p\left(
u,\epsilon \left( u\right) \right) \right\vert ^{2}+f_{\epsilon (t)}\left(
y\left( u,\epsilon \left( u\right) \right) \right) \right] du\geq
z_{\epsilon \left( 0\right) }\left( y\left( 0,\epsilon \left( 0\right)
\right) \right) .
\end{equation*}%
Moreover%
\begin{eqnarray*}
Ee^{-\alpha _{\epsilon \left( \tau ^{\ast }\right) }\tau ^{\ast
}}z_{\epsilon \left( \tau ^{\ast }\right) }\left( y^{\ast }\left( \tau
^{\ast },\epsilon \left( \tau ^{\ast }\right) \right) \right) &=&Ee^{-\alpha
_{\epsilon \left( \tau \right) }\tau }z_{\epsilon \left( \tau \right)
}\left( y\left( \tau ,\epsilon \left( \tau \right) \right) \right) \\
&=&z_{\epsilon \left( \tau ^{\ast }\right) }\left( R\right) =z_{\epsilon
\left( \tau \right) }\left( R\right) =0,
\end{eqnarray*}%
since%
\begin{equation*}
z_{i}\left( R\right) =0\text{, }i=1,2.
\end{equation*}%
This together with $y^{\ast }\left( 0,\epsilon \left( 0\right) \right)
=y\left( 0,\epsilon \left( 0\right) \right) $ finishes the proof.

\paragraph{Acknowledgments}

\end{document}